\documentclass{amsart}%
\usepackage{amsfonts}
\usepackage{amsmath}
\usepackage{amssymb}
\usepackage{graphicx}%
\providecommand{\U}[1]{\protect\rule{.1in}{.1in}}

\theoremstyle{plain}

\newtheorem{remark}{Remark}

\numberwithin{equation}{section}

\begin{document}

$$\large{\bf \text{Some aspects of quasi-pseudo principally injective modules}}$$

\qquad \qquad \qquad \quad Hemen Dutta$^1$, Azizul Hoque$^2$ and Samer M. Saeed$^3$\\

\noindent $^1$ Gauhati University, Faculty of Science, Department of Mathematics,
Guwahati-781014, India\\
E-mail: hemen\_dutta08@rediffmail.com\\

\noindent$^2$ Gauhati University, Faculty of Science, Department of Mathematics,
Guwahati-781014, India\\
E-mail: ahoque.ms@gmail.com\\

\noindent$^3$The Iraqi Ministry of Education, Directorate-General for Education, Wasit, Iraq\\
E-mail: samirseed@yahoo.com\\

$$\small ABSTRACT$$
In this paper, the notion of quasi-pseudo injectivity relative to a class of submodules, namely, quasi-pseudo principally injective has been studied. This notion is closed under direct summands. Several properties and characterizations have been given. In particular, we characterize Noetherian Rings and Dedekind Domains by quasi-pseudo principally injectivity.\\

\noindent 2010 AMS Classification: 16D10; 16D50; 16D80; 16D99.\\
Keywords: Quasi-pseudo principally injective module; Pseudo principally injective module; Noetherian ring; Dedekind domain.

\section{Introduction}

Throughout this paper, $R$ represents an associative ring with identity. It is well-known that $R$-modules are unitary right $R$-modules. Let $M$ be an $R$-module, a module $N$ is called $M$-generated, if there is an epimorphism $M^{(I)}\rightarrow N$ for some index set $I$. If $I$ is finite, then $M$ is called finitely $M$-generated. In particular, a submodule $N$ of $M$ is called an $M$-cyclic submodule of $M$, if it is isomorphic to $\frac{M}{L}$, for some submodule $L$ of $M$. Hence, any $M$-cyclic submodule $N$ of $M$ can be considered as the image of an endomorphism of $M$ \cite{W1991}. An $R$-module $M$ is called epi-retractable, if every submodule of $M$ is a homomorphic image of $M$ \cite{GV2009}. A submodule $N$ of an $R$-module $M$ is said to be a direct summand of the $R$-module $M$, if there exists a submodule $L$ of $M$ $M=N\oplus L$ \cite{G1976}. It is well-known that every direct summand submodule is $M$-cyclic. A submodule $N$ of an $R$-module $M$ is said to be essential submodule of the $R$-module $M$, if $N$ has non-zero intersection with every non-zero submodule of $M$. A non-zero $R$-module $U$ is said to be uniform if every non-zero submodule of $U$ is essential in $U$ (\cite{G1976}). We shall use $\vartheta(R)$ to stand for the set of all essential right ideals of the ring $R$. Given any $R$-module $M$, we set $Z(M)$ = $\{x\in M|xI=0$, for some $I\in \vartheta(R)\}$. An $R$-module $M$ is said to be singular if $Z(M)$ = $M$. At the other extreme, we say $M$ is non-singular if $Z(M)$ = 0 (\cite{G1976}). An $R$-module $M$ is said to be semisimple, if every submodule of $M$ is direct summand (\cite{G1976}).

Consider the following conditions for an $R$-module $M$:
\begin{itemize}
\item[(C$_1$)] Every submodule of $M$ is essential in a direct summand of $M$.
\item[(C$_2$)] If a submodule of $M$ is isomorphic to a direct summand of $M$, then it is a direct summand of $M$ itself.
\item[(C$_3$)] If $A$ and $B$ are direct summands of $M$ with $A\cap B$ = 0, then $A\oplus B$ is a direct summand of $M$.
\end{itemize}
An $R$-module $M$ is called a $CS$-module or an extending module if $M$ satisfies (C$_1$)  \cite{MM1990}. A $CS$-module $M$ which satisfies (C$_2$)-condition is called continuous (\cite{MM1990}). Moreover a $CS$-module $M$ which satisfies (C$_3$)-condition is called quasi-continuous(\cite{MM1990}).
Let $M$ and $N$ be two $R$-modules. Then $N$ is called $M$-injective, if for every submodule $A$ of $M$, any $R$-homomorphism from $A$ to $N$ can be extended to an $R$-homomorphism from $M$ to $N$. An $R$-module $N$ is called injective, if it is $M$-injective for all $R$-module $M$. An $R$-module $M$ is called quasi-injective, if it is $M$-injective. In 2005, Dinh \cite{D2005} studied the notion of pseudo-$M$-injective modules (the original terminology is $M$-pseudo-injective), following which a right $R$-module $N$ is called pseudo-$M$-injective if for every submodule $A$ of $M$, any monomorphism $\alpha : A\longrightarrow N$ can be extended to a homomorphism $\beta : M\rightarrow N$. A right $R$-module $N$ is called pseudo-injective if $N$ is pseudo $N$-injective \cite{SJ1967}. Sanh et.al. \cite{SSDW1999} studied the structure of quasi principally injective modules and presented several characterizations and properties. An $R$-module $N$ is called $M$-principally injective, if every $R$-homomorphism from an $M$-cyclic submodule $A$ of $M$ to $N$ can be extended to an $R$-homomorphism from $M$ to $N$. A module $M$ is called quasi principally (or semi) injective, if it is $M$-principally injective. In 2010, Patel et. al. \cite{PPGK2010} also studied the notion of quasi-principally-injective modules for some relevant properties. In this paper, we shall investigate some aspects of the notion of quasi-pseudo principally injectivity (in short, quasi-pp-injective module).

Let $M$ and $N$ be two right $R$- modules. $M$ is said to be pseudo principally-$N$-injective (or pseudo $N$- $p$-injective), if for each $M$-cyclic submodule $A$ of $N$, every $R$-monomorphism from $A$ to $M$ can be extended to an $R$-homomorphism from $N$ to $M$. A right $R$-module $M$ is called quasi-pseudo principally injective (or quasi-$pp$-injective), if it is pseudo principally $M$- injective (or pseudo $M$-$p$-injective). A right $R$-module $N$ is called pseudo principally- injective (or pseudo $p$-injective) if it is pseudo $R_R$-principally injective (or pseudo $R_R$- $p$-injective). A ring $R$ is called right self-pseudo principally-injective (or self-pseudo $p$- injective) if $R_R$ is quasi-$pp$-injective as a right $R$- module. For some additional details, please refer to Sanh et. al. \cite{SSDW1999}, Baupradist et. al. \cite{BHS2011}, etc.

\begin{remark} \label{rem1} We have the following immediate results:
\begin{enumerate}
  \item Every pseudo-injective module is pseudo-principally-injective (of course, these two notions coincide for epi-retractable R-modules).
  \item Every quasi-principally-injective module is pseudo-principally-injective.
  \item An R-isomorphic module to quasi-pp-injective is quasi-pp-injective.
  \item An R-isomorphic module to pseudo-principally-injective is pseudo-principally-injective.
  \item If $N_1$ and $N_2$ are isomorphic $R$-modules and $M$ is a pseudo-$p$-$N_1$-injective $R$-module, then $M$ is a pseudo-$p$-$N_2$-injective $R$-module. 
\end{enumerate}
\end{remark}

\section{\bf Main Results and Discussion}

In this section, we shall discuss and prove our main results. 

\noindent \textbf{Definition 2.1.}\ \
An $R$-module $M$ is called a multiplication module if each submodule $N$ of $M$ is of the form $MI$ for some right ideal $I$ of $R$ \cite{B1981}.

\noindent \textbf{Proposition 2.2.}\ \
Let $M$ be a multiplication quasi-pp-injective module. Then every $M$-cyclic submodule of $M$ is quasi-pp-injective.

\noindent \textbf{Proof:}\ \
Let $N$ be an $M$-cyclic submodule of $M$. To show $N$ is quasi-pp-injective.
Let us take $L$ be an $N$-cyclic submodule of $N$. Then $L$ is also $M$-cyclic submodule of $M$.
Consider a monomorphism $\alpha:L\longrightarrow N$. Then $\alpha$ is a monomorphism from $L$ into $M$. Thus by the quasi-pp-injectivity of $M$, there exists an $R$-homomorphism $\beta: M\longrightarrow M$ such that $\beta_{|_L}=\alpha$.

since $M$ is multiplicative, $N$ can be written as $N=MI$ for some ideal of $R$. Thus $\beta(N)
=\beta(MI)=\beta(M)I\subseteq MI=N$. Therefore $\beta^*=\beta_{|N}:N\longrightarrow N$ is a an $R$-homomorphism. Moreover $\beta^*_{|_L}=\alpha$. This completes the proof.
\\
In the following result, we show that the connection between quasi-pp-injectivity and semi-simplicity for self-similar modules. 

\noindent \textbf{Definition 2.3.}\ 
An $R$-module is said to be self-similar if every submodule of $M$ is isomorphic to $M$ \cite{RS2009}.

\noindent \textbf{Proposition 2.4.}\ 
Let $M$ be self-similar $R$-module. Then $M$ is quasi-pp-injective if and only if $M$ is semi-simple.
\\
As a consequence we derive the following result.

\noindent \textbf{Corollary 2.5.}\
If $R_R$ is self-similar, then $R_R$ is quasi-pp-injective if and only if $R_R$ is semi-simple. 
\\
Using this results we can conclude that if $R_R$ is self-similar and semi-simple, then the ring $R$ is self-pseudo principally injective (self-pseudo-p-injective).

\noindent \textbf{Proposition 2.6} \cite{CPG2009}.\ 
Every quasi-pp-injective module satisfies $C_2$.
\\
Quasi-pp-injectivity is not closed under direct sum, as we see in the following:\\
$R= \begin{bmatrix}
F & F\\
F & 0
\end{bmatrix}$,
$A=\begin{bmatrix}
0 & 0 \\
F & 0
\end{bmatrix}$,
$B=\begin{bmatrix}
F & F \\
0 & 0
\end{bmatrix}$ and
$C=\begin{bmatrix}
0 & F \\
0 & 0
\end{bmatrix}$, where $F=\mathbb{Z} / 2 \mathbb{Z}$.
It is easy to see that the right $R$-modules $A$ and $B$ are quasi-pp-injective (in fact they are quasi-injective). However $R=A\oplus B$ is not quasi-pp-injective, since otherwise $R$ satisfies (C$_2$) by proposition 2.6., but $A$ is isomorphic to $C$ and $C$ is not a direct sum in $R$.
\\
However in case direct summand we have the following result.

\noindent \textbf{Proposition 2.7.}\ 
Every direct summand of a quasi-pp-injective (respectively pseudo-p-injective) $R$-module is quasi-pp-injective (respectively pseudo-p-injective).

\noindent \textbf{Theorem 2.8.}\ 
The following statements are equivalent for quasi-pp-injective $R$-module $M$.
\begin{enumerate}
\item $M$ is continuous.
\item $M$ is quasi-continuous.
\end{enumerate}
We now explore more properties of quasi-pp-injective modules.

\noindent \textbf{Proposition 2.9.}\ 
Let $R$-module $M$ be direct sum of two submodules $M_1$ and $M_2$. If $M$ is quasi-pp-injective modules then $M_1$ and $M_2$ are mutually-principally injective.

\noindent \textbf{Proof:}\ 
Let $A$ be an $M_2$-cyclic submodule of $M_2$ and $f:A\longrightarrow M_1$ be an $R$-monomorphism. Define $g:A\longrightarrow M=M_1 \oplus M_2$ by $g(a)=(f(a), a), \forall a\in A$. Then $g$ is an $R$-monomorphism. By the hypothesis $M$ is pseudo $M_2$-p-injective and thus $g$ extends to an endomorphism $h$ of $M$. Put $h^*=h_{|M_2}$. Let $\pi : M\longrightarrow M_1$ be the natural projection of $M$ onto $M_1$. Then $\pi\circ h^*:M_2 \longrightarrow M_1$ extends $f$. This shows that $M_1$ is principally $M_2$-injective module. Similarly we can show that $M_2$ is principally $M_1$-injective module.  
\\
Recently Abbas and Saied \cite{AS2012} introduced the concept of IC-pseudo injective modules. The next result characterizes quasi-pp-injective modules in terms of IC-pseudo injective modules.
 
\noindent \textbf{Definition 2.10.} \ 
Let $M$ and $N$ be two $R$-modules. If for every submodule $A$ which is isomorphic to closed submodule of $M$, any $R$-monomorphism from $A$ into $N$ can be extended to an $R$-homomorphism from $M$ into $N$, then $N$ called IC-pseudo $M$-injective. An $R$-module $M$ is called IC-pseudo-injective if it is ic-pseudo-M-injective. A submodule $N$ of $M$ is called IC-submodules if $N$ is isomorphic to a closed submodule of $M$.

\noindent \textbf{Proposition 2.11.}\ 
Every uniform non-singular IC-pseudo injective module is quasi-pp-injective.

\noindent \textbf{Proof:}\ 
Let $M$ be a uniform non-singular IC-pseudo injective module. Let $N$ be an $M$-cyclic submodule of $M$. and let $\psi: N\longrightarrow M$ be an $R$-monomorphism. By hypothesis, there is an epimorphism $\phi: M\longrightarrow N$ and $\frac{M}{ker f}$ is non-singular, and thus $ker =0$. Therefore $N$ is isomorphic to $M$ and thus $N$ is IC-submodule of $M$. Hence there exists an $R$-homomorphism $\psi^*: M\longrightarrow M$ such that $\psi^*_{|_N}=\psi$.
\\
It is well known that an $R$-module $M$ is injective if and only if $M$ is $N$-injective for each $R$-module $N$.
 
\noindent \textbf{Proposition 2.12.}\ 
The following statements are equivalent for an $R$-module M.
\begin{enumerate}
\item $M$ is injective.
\item $M$ is pseudo principally $N$-injective, for each $R$-module $N$.
\end{enumerate}
\noindent \textbf{Proof:}\ 
$(1)\Longrightarrow (2)$ is obvious.
\\
$(2)\Longrightarrow (1)$: Let $E=E(M)$ be the injective hull of $M$. Then $M$ is $M$-cyclic submodule of $M \oplus E$. For the inclusion mapping $i: M\longrightarrow E$ and the natural injection $j:E\longrightarrow M\oplus E$, the pseudo-principally-$M\oplus E$-injectivity of $M$ implies that the identity $R$-homomorphism of $M$ can be extended to an $R$-homomorphism $f:M\oplus E \longrightarrow M$. This shows that $M$ is a direct summand of $E$ and hence $M$ is injective.
\\
In \cite{JLP1995}, Jain and L\'{o}pez-Permouth recall that an $R$-module $M$ is weakly injective if for each finitely generated $R$-module $N$ and for each $R$-homomorphism $\phi: N \longrightarrow E(M)$, there exists a submodule $X$ of $E(M)$ such that $\phi(N)\subset X \cong M$. Alternatively, $M$ is weakly injective if for each finitely generated $R$-module $N\ll E(M)$, we have $N\ll X\ll E(M)$ for some $X\cong M$. 

\noindent \textbf{Proposition 2.13.}\ 
Let $M$ be a finitely generated right $R$-module. Then $M$ is injective if and only if it is weakly injective and quasi-pp-injective.

\noindent \textbf{Proof:}\ 
We need only to prove the sufficiency.\\
Suppose $M$ is weakly injective and quasi-pp-injective. Let $m\in E(M)$. By weakly injectivity of $M$, there exists $N\ll E(M)$ such that $M+mR\ll N\cong M$. Also by quasi-pp-injectivity of $M$ and remark 1 (statement (3)), $N$ is quasi-pp-injective. By Proposition 2.8, $M$ is a direct summand of $N$, that is, there exists $K\ll N$ such that $N=M\oplus K$. As $M$ is an essential submodule of $N$, $N=M$ and thus $m\in M$. Therefore $M=E(M)$ which implies that $M$ is injective.
\\
Han and Chol \cite{HC1995} established that every pseudo principally injective module is a direct injective and every direct injective module is divisible. In this connection we have the following result which is also hold for quasi-pp-injective modules.

\noindent \textbf{Proposition 2.14.}\ 
Every pseudo-principally injective module is divisible.
\\
In the following we show that the distinction between quasi-pp-injectivity and divisibility vanishes over Dedekind domain. Rotmann \cite{R2000}(Theorem 4.24) recall that a domain $R$ is Dedekind ring if every divisible $R$-module is injective.   
\\
We now provide a characterization of Dedekind rings in terms of quasi-pp-injective modules.

\noindent \textbf{Proposition 2.15.}\ 
The following statements are equivalent for a ring $R$.
\begin{enumerate}
\item $R$ is Dedekind domain;
\item Every divisible $R$-module is quasi-pp-injective.
\end{enumerate}

\noindent \textbf{Proof:}\
$(1)\Longrightarrow (2)$ holds by Theorem 4.24 of \cite{R2000}.\\
$(2)\Longrightarrow (1)$: Let $M$ be a divisible module and $E$ be the injective hull of $M$. Then $E$ is divisible and hence $M\oplus E$ is also divisible. Therefore by Proposition 2.7 and Proposition 2.12, $M$ is injective.
\\
We now provide a characterization of Noetherian rings in terms of quasi-pp-injective modules. 
We recall that an $R$-module $M$ is $F-$injective if for every finitely generated ideal of $I$ of $R$, every $R$-homomorphism of $I$ into $M$ can be extended to an $R$-homomorphism of from $R$ into $M$ \cite{M1989}.
 
\noindent \textbf{Proposition 2.16.}\ 
The following conditions are equivalent for a ring $R$.
\begin{enumerate}
\item $R$ is Noetherian ring;
\item Every $F$-injective $R$-module is injective;
\item Every $F$-injective $R$-module is quasi-pp-injective.
\end{enumerate}
\noindent \textbf{Proof:}\
$(1) \Longrightarrow (2)$ and $(2) \Longrightarrow (3)$ are obvious.
\\
$(3) \Longrightarrow (2)$:
Let $M$ be an $F$-injective $R$-module and $E$ be the injective hull of $M$. We write $Q=M\oplus E$ and $S=\lbrace (m, 0): m\in M\rbrace$. Then $S$ is a direct summand of $Q$ and $S\cong M$. Since $Q$ is a direct sum of two $F-$injective $R$-modules, $Q$ is $F$-injective. By hypothesis, $Q$ is quasi-pp-injective. \\
Let $i: M \longrightarrow E$ be the inclusion map and $k: E\longrightarrow Q$ be the canonical injection. Then $\phi=k\circ i: M\longrightarrow Q$ is an $R$-monomorphism. By the quasi-pp-injectivity of $Q$, for the canonical injection $j:M\longrightarrow Q$, there exists an $R$-homomorphism $h:Q\longrightarrow Q$ such that $h\circ \phi=j$. If $\pi:Q\longrightarrow M$ is the canonical projection, then $\psi=\pi\circ h\circ k: E\longrightarrow M$ such that $\psi\circ i=\pi\circ j$ is the identity map on $M$. Therefore $M$ is a direct summand of $E$ which yeilds $M$ is injective.  
\\
$(2) \Longrightarrow (1)$: Assume (2). Since any direct sum of $F$-injective $R$-modules is $F$-injective, by (2) any direct sum of injective $R$-modules is injective. Therefore by (\cite{R2000}, P. 82), $R$ is Noetherian.


\begin{thebibliography}{99}

\bibitem{AS2012} {\sc M.S. Abbas and S.M. Saied}, {\em IC-pseudo-injective modules}, International Journal of Algebra, \textbf{6(6)} (2012), 255-264.

\bibitem{B1981} {\sc A. Barnard}, {\em Multiplication modules}, J. Algebra, \textbf{71} (1981), 174-178.

\bibitem{HC1995} {\sc Chang-woo Han and Su-Jeong Chol}, {\em Generalizations of the quasi-injective modules}, Comm. Korean Math. Soc., \textbf{10(4)} (1995), 811-813.

\bibitem{CPG2009} {\sc  A.K. Chaturvedi, B.M. Pandeya and A.J. Gupta}, {\em Quasi-pseudo-principally injective modules}, Algebra Colloq., \textbf{16(3)} (2009), 397- 402.

\bibitem{D2005} {\sc H.Q. Dinh}, {\em A note on pseudo-injective modules}, Comm. Algebra, \textbf{33} (2005), 361-369.

\bibitem{GV2009} {\sc A. Ghorbani and M.R. Vedadi}, {\em Epi-retractable modules and some applications}, Bul. Iranian Math. Soc., \textbf{35(1)} (2009), 155-166.

\bibitem{G1976} {\sc K.R. Goodearl}, {\em Ring Theory, Nonsingular Rings and Modules}, Marcel Dekker. Inc., New York, 1976..

\bibitem{JLP1995} {\sc S.K. Jain and S.R. L\'{o}pez-Permouth}, {\em Weakly-injective modules over  hereditary noetherian prime rings}, J. Austral. Math. Soc.(Series A), \textbf{58} (1995), 287-297.

\bibitem{MM1990} {\sc S.H. Mohamed and B.J. M\"{u}ller}, {\em Continuous and Discrete Modules}, London Mathematical Society Lecture Note Series 14, Cambridge University Press, 1990.

\bibitem{M1989} {\sc R.Y.C. Ming}, {\em On regular rings and self-injective rings, IV}, Publications De l'institut Mathematique, Nouvelle serie tome, \textbf{45(59)} (1989), 65-72.

\bibitem{RS2009} {\sc V.S. Rodrigues and A.A. Sant'Ana}, {\em A note on a problem due to Zelmanowitz}, Algebra and Discrete Mathematics, \textbf{3} (2009), 85-93.

\bibitem{R2000} {\sc J.J. Rotman}, {\em An Introduction to Homological Algebra}, Academic Press, New York, 2000.

\bibitem{SSDW1999} {\sc N. V. Sanh, K. P. Shum, S. Dhompongsa and S. Wongwai}, {\em On Quasi Principally Injective Modules}, Algebra Colloq., \textbf{6(3)} (1999), 269-276.

\bibitem{SJ1967} {\sc S. Singh and S. K. Jain}, {\em On Pseudo-injective modules and self pseudo-injective rings}, J. Math. Sci., \textbf{2} (1967), 23-31.

\bibitem{BHS2011} {\sc S. Baupradist, H.D. Hai and N.V. Sanh}, {\em On pseudo-p-injectivity}, Southeast Asian Bull. Math., \textbf{35} (2011), 21-27.

\bibitem{PPGK2010} {\sc M.K. Patel, B.M. Pandeya, A.J. Gupta and V. Kumar}, {\em Quasi principally injective modules}, International Journal of Algebra, \textbf{4(26)} (2010), 1255 - 1259.

\bibitem{W1991} {\sc R. Wisbauer}, {\em Foundations of Module and Ring Theory}, Revised and translated from the 1988 German edition, Algebra, Logic and Applications, vol. 3, Gordon and Breach Science Publishers, Philadelphia, PA, 1991.

\end{thebibliography}
\end{document}